\documentclass[12pt]{amsart}
\usepackage{amssymb}
\theoremstyle{plain}
 \newtheorem{thm}{Theorem}[section]
 
 \newtheorem{lem}[thm]{Lemma}

\theoremstyle{definition}

\theoremstyle{remark}

\begin{document}
\title[A characterization of 
the operator-valued triangle 
equality]
{A characterization of 
the operator-valued triangle 
equality}
\author[Tsuyoshi Ando and 
Tomohiro Hayashi]{{
Tsuyoshi Ando and  Tomohiro Hayashi} }
\address[Tsuyoshi Ando]{Hokkaido University (Emeritus)}
\address[Tomohiro Hayashi]
{Graduate School of Mathematics, Kyushu University, 
33, Fukuoka, 812-8581, Japan}
\email[Tsuyoshi Ando]{ando@es.hokudai.ac.jp}
\email[Tomohiro Hayashi]{hayashi@math.kyushu-u.ac.jp}
\baselineskip=17pt

\maketitle

\begin{abstract}
We will show that 
for any two bounded linear operators 
$X,Y$ on a 
Hilbert space ${\frak H}$, if they satisfy the triangle 
equality $|X+Y|=|X|+|Y|$, there 
exists a partial isometry $U$ on ${\frak H}$ 
such that $X=U|X|$ and $Y=U|Y|$. 
This is a generalization of Thompson's 
theorem 
to the matrix case proved by using 
a trace. 
\end{abstract}



\section{Introduction}
In his interesting paper~\cite{T1,T2}, Thompson showed the triangle 
inequality for two matrices: 
for any two matrices 
$X$ and $Y$, there exist unitaries 
$V$ and $W$ such that 
$|X+Y|\leq V|X|V^{*}+W|Y|W^{*}$, 
where $|X|=(X^{*}X)^{1/2}$. 
(See~\cite{AAP} for a generalization of the theorem 
to bounded linear operators on infinite dimensional 
Hilbert spaces.) 
Moreover 
in~\cite{T3} he characterized the condition 
when the triangle equality holds: he showed that 
for two matrices $X$ and $Y$, 
there exists two unitary matrices 
$V,W$ such that 
the equality 
$|X+Y|=V|X|V^{*}+W|Y|W^{*}$ holds 
if and only if there is a 
unitary matrix $U$ such that 
$X=U|X|$ and $Y=U|Y|$.

In the proof of Thompson's theorem, 
the (finite) trace plays a crucial role. 
Thus his argument does not work in the 
infinite dimensional case. The aim of this paper 
is to overcome this difficulty. 
We shall generalize this theorem 
to operators on 
infinite-dimensional Hilbert spaces. 

It is easy to see that Thompson's theorem 
fails in the infinite dimensional setting. 
Indeed we can easily find two operators 
$X$, $Y$ and unitaries $V$, $W$ 
satisfying $|X+Y|=V|X|V^{*}+W|Y|W^{*}$, 
but there is no partial isometry 
$U$ such that 
$X=U|X|$ and $Y=U|Y|$. 
For example, consider an orthogonal 
projection $P$ and its subprojection 
$Q$. Suppose that all projections 
$P$, $Q$ and $P-Q$ are 
infinite rank. Then there are 
(at least countably many) 
unitaries $V,W$ such that 
the equality 
$|P+(-Q)|=V|P|V^{*}+W|(-Q)|W^{*}$ 
holds. In this case there is no 
partial isometry $U$ satisfying 
$P=U|P|$ and $(-Q)=U|(-Q)|$. 

To avoid such cases, we consider the following 
problem: 
for two bounded linear operators 
$X$ and $Y$ on an infinite dimensional 
Hilbert space ${\frak H}$, 
if the equality 
$|X+Y|=|X|+|Y|$ holds, 
does there exists a 
partial isometry $U$ 
on ${\frak H}$ 
such that 
$X=U|X|$ and $Y=U|Y|$? 
(Here we remark even in this case 
we cannot choose $U$ as a unitary. 
For example, consider the unilateral 
shift $S$ 
(\cite{H}, p40). Then obviously the 
equality 
$|S+S|=|S|+|S|$ holds 
but there exists no 
unitary $U$ satisfying 
$S=U|S|$.) 
The main result of this paper is to 
answer this problem affirmatively.

We would like to thank Professor 
Yoshimichi Ueda who suggested this problem 
to us. We also would like to 
express their hearty gratitude to 
Professors Fumio Hiai 
and Atsushi Uchiyama for valuable 
comments. 

The second author was supported by 
Grant-in-Aid for Young Scientists 
(B) of Japan Society for 
the Promotion of Science.

\section{Main Result}

Throughout this paper we assume that 
the readers are familiar with 
basic notations and results on the 
Hilbert space operator 
theory. We refer to 
Halmos's book~\cite{H}. 

We denote 
by ${\frak H}$ an infinite dimensional 
complex Hilbert space. 
The inner product 
on ${\frak H}$ is denoted by 
$\langle x,y\rangle$ 
for two vectors $x,y\in {\frak H}$. 
Thus the Hilbert space norm $||\cdot||$ 
is defined by 
$||x||=\langle x,x\rangle^{1/2}$. 
We denote by $B({\frak H})$ 
the set of 
all bounded linear operators on ${\frak H}$. 
The operator norm on $B({\frak H})$ 
is denoted by $||\cdot||$. 
The identity 
operator on ${\frak H}$ is denoted by $I$. 
An operator $W\in B({\frak H})$ is 
called a contraction 
if $||W||\leq 1$. 
An operator $A\in B({\frak H})$ is said 
to be positive 
if $\langle Ax,x\rangle\geq 0$ for any 
$x\in {\frak H}$. 
We denote $A\geq 0$ if $A$ is positive. 
For two selfadjoint operators $A,B\in B({\frak H})$, 
we denote $A\leq B$ if $B-A\geq 0$. 
A positive operator 
$A\in B({\frak H})$ has the unique square 
root $A^{1/2}\in B({\frak H})$ which satisfies 
$(A^{1/2})^{2}=A$ and $A^{1/2}\geq 0$. 
For $X\in B({\frak H})$, we define its 
``absolute value'' by $|X|=(X^{*}X)^{1/2}$. 
An operator $U\in B({\frak H})$ is called 
a partial isometry if 
$U^{*}U$ is an orthogonal projection. 
For $X\in B({\frak H})$, we can always express 
$X$ as a product $X=U|X|$, where $U\in B({\frak H})$ is a 
partial isometry. Then $U$ satisfies that 
$U^{*}U|X|=|X|$. This expression is called 
a polar decomposition of $X$. 
We define the kernel ${\rm ker}(X)$ 
and the range ${\rm ran}(X)$ of $X$ by 
${\rm ker}(X)=\{a\in {\frak H}:\ 
Xa=0\}$ and 
${\rm ran}(X)=\{Xa:\ 
a\in {\frak H}\}$ respectively. 
We denote the set of all 
real numbers by ${\Bbb R}$ 
and the set of all 
complex numbers by ${\Bbb C}$. 
For $z\in {\Bbb C}$ its real 
part is denoted by 
${\rm Re}(z)$. 
The spectrum of $X\in B({\frak H})$ 
is defined by 
$\sigma(X)=\{z\in {\Bbb C}:\ 
X-zI\ 
{\text{does not have an inverse in}}\ B({\frak H})\}$. 
The numerical range of $X\in B({\frak H})$ 
is also defined by 
${\Bbb W}(X)=
\{\langle Xa,a\rangle :\ a\in {\frak H},\ 
||a||=1 \}$.

\begin{lem}
Let $A,B\in B({\frak H})$ be two positive 
operators and let $V,W\in B({\frak H})$ be 
partial isometries 
such that 
$V^{*}VA=A$ and $W^{*}WB=B$. 
Then there exists a partial isometry 
$U\in B({\frak H})$ such that 
$UA=VA$ and $UB=WB$ 
if and only if 
the equality 
$A(V^{*}W-I)B=0$ holds.
\end{lem}

\begin{proof}
The ``only if'' part is obvious. 

Assume that the equality 
$A(V^{*}W-I)B=0$ holds. 
Since for any $x,y\in {\frak H}$ 
\begin{align*}
||VAx+WBy||^{2}&=||Ax||^{2}
+||By||^{2}+2{\rm{Re}}
\langle AV^{*}WBy,x\rangle\\
&=||Ax||^{2}
+||Bx||^{2}+2{\rm{Re}}
\langle ABy,x\rangle\\
&=||Ax+By||^{2}, 
\end{align*}
we can define the desired partial isometry 
$U$ by $U(Ax+By)=VAx+WBy$. 
\end{proof}

\begin{lem}
Let $W\in B({\frak H})$ be a contraction. 
Then for any 
operator $D\in B({\frak H})$ we have 
$$\sigma(D^{*}(W-I)D)\subset 
\{z\in {\Bbb C}:\ |z+||D||^{2}|\leq ||D||^{2}\},$$ 
and hence 
$$\sigma(D^{*}(W-I)D)\cap i{\Bbb R}\subset \{0\}.$$ 
\end{lem}

\begin{proof}
For any unit vector $x\in {\frak H}$, 
we see that 
$$
\langle D^{*}(W-I)Dx,x\rangle
=||Dx||^{2}\langle (W-I)\frac{Dx}{||Dx||},\frac{Dx}{||Dx||}\rangle.
$$ 
Since $W$ is a contraction, 
the numerical range ${\Bbb W}(W-I)$ is 
contained in the set 
$\{z\in {\Bbb C}:\ |z+1|\leq 1\}$. 
Combining this with $||Dx||^{2}\leq ||D||^{2}$, 
we conclude that 
$\langle D^{*}(W-I)Dx,x\rangle 
\in \{z\in {\Bbb C}:\ |z+||D||^{2}|\leq ||D||^{2}\}$. 
In other words, the numerical range 
${\Bbb W}(D^{*}(W-I)D)$ is contained in 
$\{z\in {\Bbb C}:\ |z+||D||^{2}|\leq ||D||^{2}\}$. 

Recall that in general, for any operator 
$X\in B({\frak H})$ the spectrum 
$\sigma(X)$ is a subset of the 
closure of its numerical range 
${\Bbb W}(X)$. 
(See~\cite{H}, p.111, Problem 169.) 
Hence we get 
$$\sigma(D^{*}(W-I)D)\subset 
\{z\in {\Bbb C}:\ |z+||D||^{2}|\leq ||D||^{2}\}.$$
\end{proof}

\begin{thm} 
For two bounded linear operators $X,Y\in B({\frak H})$, 
the triangle equality 
$|X+Y|=|X|+|Y|$ holds 
if and only if there exists a 
partial isometry $U$ such that 
$X=U|X|$ and $Y=U|Y|$. 
\end{thm}

\begin{proof} 

For notational simplicity 
we write 
$A\equiv |X|$ and $B\equiv |Y|$. 

First assume that there exists a 
partial isometry $U$ such that 
$X=UA$ and $Y=UB$. Then since 
$$U^{*}UA=A\ {\text{and}}\ U^{*}UB=B,$$ 
we have 
\begin{align*}
|X+Y|^{2}&=(AU^{*}+BU^{*})(UA+UB)\\
&=A^{2}+AB+BA+B^{2}=(A+B)^{2}\\
&=(|X|+|Y|)^{2}
\end{align*}
and so $|X+Y|=|X|+|Y|$. 

Assume conversely 
that the triangle equality 
$|X+Y|=|X|+|Y|$ holds. 
Let $X=V|X|=VA$ and $Y=W|Y|=WB$ be 
the polar decompositions 
of $X$ and $Y$ respectively. 
Then the triangle equality implies 
$$(AV^{*}+BW^{*})(VA+WB)=(A+B)^{2}.$$ 
Combining this with the relations $V^{*}VA=A$ 
and $W^{*}WB=B$ we have 
$$A(V^{*}W-I)B+B(W^{*}V-I)A=0.$$ 
This means that the operator 
$iA(V^{*}W-I)B$ is selfadjoint. 
We claim 
$$iA(V^{*}W-I)B=0. \eqno{(\dagger)}$$ 

Since $A,B\leq A+B$, there are 
uniquely contractions $K,L$ 
such that 
$$A^{1/2}=K(A+B)^{1/2}
\ {\text{and}}\ 
B^{1/2}=L(A+B)^{1/2}
\eqno{(1)}$$ 
and ${\rm ker}(A+B)^{1/2}
\subset 
{\rm ker}(A)\cap{\rm ker}(B)$. 
By taking adjoints of the 
both sides of (1) we also have 
$$A^{1/2}=(A+B)^{1/2}K^{*}
\ {\text{and}}\ 
B^{1/2}=(A+B)^{1/2}L^{*}.
\eqno{(2)}$$ 
Let $P$ be the orthogonal projection 
to the closure of 
${\rm ran}(A+B)^{1/2}$. 
(Remark that $I-P$ 
is the orthogonal projection 
to ${\rm ker}(A+B)^{1/2}$ 
by selfadjointness of 
$A+B$. ) 
Since by (1) and (2) 
$$PK^{*}K=K^{*}KP=K^{*}K
\ {\text{and}}\ 
PL^{*}L=L^{*}LP=L^{*}L$$ 
and 
$$(A+B)^{1/2}(K^{*}K+L^{*}L-I)(A+B)^{1/2}=0$$ 
we conclude 
$$K^{*}K+L^{*}L=P,$$ 
because ${\rm ran}(A+B)^{1/2}$ is 
dense in ${\rm ran}(P)$. 
This implies that 
$K^{*}K$ and $L^{*}L$ commute, 
and hence $(L^{*}L)(K^{*}K)$ is a 
positive
operator. 

Now by (1) and (2) we can write 
\begin{align*}
iA&(V^{*}W-I)B\\
&=i(A+B)^{1/2}K^{*}K(A+B)^{1/2}(V^{*}W-I)
(A+B)^{1/2}L^{*}L(A+B)^{1/2}.
\end{align*}
Since $iA(V^{*}W-I)B$ is selfadjoint, 
as in the proof of $K^{*}K+L^{*}L=P$ 
we see that 
$iK^{*}K(A+B)^{1/2}(V^{*}W-I)
(A+B)^{1/2}L^{*}L$ 
is also selfadjoint, and for
$(\dagger)$ it suffices to show that 
$$i(K^{*}K)(A+B)^{1/2}(V^{*}W-I)
(A+B)^{1/2}(L^{*}L)=0.\eqno{(\ddagger)}$$ 
Then since in general 
a selfadjoint operator whose 
spectrum is a one-point set 
$\{0\}$ must be $0$ (See~\cite{H}, p.61.), 
for $(\ddagger)$ it suffices to show 
$$\sigma(i(K^{*}K)(A+B)^{1/2}(V^{*}W-I)
(A+B)^{1/2}(L^{*}L))=\{0\}.$$ 
Since every spectre of a selfadjoint operator is real, 
it is further reduced to showing that 
$$\sigma((K^{*}K)(A+B)^{1/2}(V^{*}W-I)
(A+B)^{1/2}(L^{*}L))\cap (i{\Bbb R})= \{0\}. 
\eqno{(\sharp)}$$  
To prove $(\sharp)$ we use a general fact that 
for two bounded linear operators $S,T$ 
$$\sigma(ST)\setminus\{0\}
=\sigma(TS)\setminus\{0\}. \eqno{(3)}$$ 
(See~\cite{H}, p.39, Problem 61.) 
Then we have by (3) 
\begin{align*}
\sigma&((K^{*}K)(A+B)^{1/2}(V^{*}W-I)
(A+B)^{1/2}(L^{*}L))\setminus\{0\}\\
&=\sigma((V^{*}W-I)
(A+B)^{1/2}(L^{*}L)(K^{*}K)(A+B)^{1/2})\setminus\{0\}.
\end{align*}
Since $(L^{*}L)(K^{*}K)$ is a 
positive operator, so is 
$(A+B)^{1/2}(L^{*}L)(K^{*}K)(A+B)^{1/2}$. 
Let 
$$D\equiv [(A+B)^{1/2}(L^{*}L)(K^{*}K)(A+B)^{1/2}]^{1/2}$$ 
and apply (3) once again to get 
\begin{align*}
\sigma&((V^{*}W-I)
(A+B)^{1/2}(L^{*}L)(K^{*}K)(A+B)^{1/2})\setminus\{0\}\\
&=\sigma(D(V^{*}W-I)D)\setminus\{0\}.
\end{align*}
Therefore we have 
\begin{align*}
\sigma&((K^{*}K)(A+B)^{1/2}(V^{*}W-I)
(A+B)^{1/2}(L^{*}L))
\cap (i{\Bbb R})\setminus\{0\}\\
&= 
\sigma(D(V^{*}W-I)D)\cap (i{\Bbb R})\setminus\{0\}. 
\end{align*}
On the other hand 
since $V^{*}W$ is a contraction we have 
by Lemma 2.2 
$$\sigma(D(V^{*}W-I)D)\cap (i{\Bbb R})\subset \{0\}$$ 
which establishes 
$(\sharp)$ and hence $(\dagger)$. 

Now by Lemma 2.1 it follows 
from $(\dagger)$ that 
there is a partial isometry $U$ 
such that 
$$U|X|=UA=VA=V|X|\ 
{\text{and}}\ 
U|Y|=UB=WB=W|Y|.$$ 
This completes the proof.

\end{proof}


\begin{thebibliography}{99}
\bibitem{AAP} 
C.~A.~Akemann, J.~Anderson and G.~K.~Pedersen, 
{\it Triangle inequalities in operator algebras.} 
Linear and Multilinear Algebra 
{\bf 11} (1982), no. 2, 167--178.



\bibitem{H} 
P.~R.~Halmos, {\it A Hilbert space problem book.} 
D. Van Nostrand Co., Inc., Princeton, N.J.-Toronto, 
Ont.-London 1967 



\bibitem{T1} 
R.~C.~Thompson, {\it Convex and concave functions of 
singular values of matrix sums.} 
Pacific J. Math. 
{\bf 66} (1976), no. 1, 285--290.

\bibitem{T2} 
\underline{\phantom{aaaaa}}, 
{\it Matrix type metric inequalities.} 
Linear and Multilinear Algebra 
{\bf 5} (1977/78), no. 4, 303--319.

\bibitem{T3} 
\underline{\phantom{aaaaa}}, 
{\it The case of equality in the matrix-valued triangle inequality.} 
Pacific J. Math. 
{\bf 82} (1979), no. 1, 279--280.
\end{thebibliography}
\end{document}